\documentclass[11pt,a4paper]{article}
\usepackage[T1]{fontenc}
\usepackage[utf8]{inputenc}
\usepackage{lmodern}
\usepackage{tikz,graphicx}
\usepackage[margin=27mm]{geometry}
\usepackage{amsmath,amssymb,amsthm,mathtools}
\usepackage{authblk}
\usepackage{enumitem,microtype}
\usepackage{booktabs}
\usepackage[hidelinks]{hyperref}

\newcommand{\gd}{\dot\gamma}
\newtheorem{theorem}{Theorem}[section]
\newtheorem{lemma}[theorem]{Lemma}
\newtheorem{conjecture}[theorem]{Conjecture}
\newtheorem{corollary}[theorem]{Corollary}
\newtheorem{proposition}[theorem]{Proposition}
\newtheorem{claim}{Claim}
\newcounter{finishstep}
\setlist[enumerate]{leftmargin=*,itemsep=1pt,topsep=3pt}
\newcommand{\leafpart}[1]{Lemma~\ref{lem2}\textup{(\ref{#1})}}
\begin{document}
	\title{\bfseries A one-third bound for domination in cubic graphs of girth at least six\thanks{Leilei Zhang's work was supported by the National Natural Science Foundation of China (Grant Nos.~12601664 and 12571365), JSPS KAKENHI (Grant No.~25KF0036), the Hubei Provincial Natural Science Foundation of China (Grant No.~2025AFB309), the China Postdoctoral Science Foundation (Grant No.~2025M773113), and the Fundamental Research Funds for the Central Universities, Central China Normal University (Grant No.~CCNU24XJ026).}}
	
	\author[1]{\bfseries Kun Cheng\thanks{Email: chengkunmath@163.com}}
	\author[2]{\bfseries Changlong Shen\thanks{Email: clshen2019@foxmail.com}}
	\author[3]{\bfseries Leilei Zhang\thanks{Email: leilei-zhang-gt@outlook.com}}
	
	\affil[1]{\small Department of Mathematics, East China Normal University, Shanghai, China}
	\affil[2]{\small School of Mathematics and Statistics, Changsha University of Science and Technology, Changsha, China}
	\affil[3]{\small Faculty of Environment and Information Sciences, Yokohama National University, Yokohama, Japan}
	\date{}
	\maketitle
	\vspace{-2.71828em}
	
	\begin{abstract}
		A set $S\subseteq V(G)$ is a dominating set of a graph $G$ if every vertex outside $S$ has a neighbor in $S$. The domination number $\gamma(G)$ is the minimum size of a dominating set of $G$. In 2010, Verstra\"{e}te conjectured that every cubic graph $G$ of order $n$ and girth at least six satisfies $\gamma(G)\le n/3$. More recently, Dorbec and Henning established this bound under the additional assumption that $G$ contains neither a $7$-cycle nor an $8$-cycle. We prove Verstra\"{e}te's conjecture in full. As a corollary, we settle the $C_4$-free case of Kostochka's conjecture, which asserts that the same bound holds for all cubic bipartite graphs.
	\end{abstract}
	
	
	\section{Introduction}\label{sec1}
	
	All graphs are finite and simple unless explicitly called multigraphs. For graph-theoretic terminology not explained in this paper, we refer the reader to~\cite{BM}. A graph is \emph{cubic} if every vertex has degree three. The \emph{girth} of a graph is the length of a shortest cycle. For a forest, we consider its girth to be infinite. For $v\in V(G)$, let $N(v)$ denote its neighborhood. For $A\subseteq V(G)$, let $N(A)=(\bigcup_{v\in A}N(v))\setminus A$ and let $G-A$ denote the graph obtained by deleting $A$. We use $\delta(G)$ for the minimum degree of $G$. The length of a path or cycle is its number of edges.  
	
	A major open problem in domination theory is to determine the maximum domination number among graphs on $n$ vertices. In 1996, Reed~\cite{Reed} proved the following classical result.
	\begin{theorem}[Reed~\cite{Reed}]
		Let $G$ be a graph of order $n$ with minimum degree at least three.
		Then $\gamma(G)\le 3n/8$.
	\end{theorem}
	For connected cubic graphs of order $n\ge 10$, Kostochka and Stodolsky~\cite{KS09} obtained the stronger bound $\gamma(G)\le 4n/11$. Kostochka and Stocker~\cite{KStocker} later improved it to $\gamma(G)\le \lfloor 5n/14\rfloor$, which remains the best known general upper bound for this class. In 1996, Reed~\cite{Reed} conjectured that $\gamma(G)\le \lceil n/3\rceil$ for every connected cubic graph of order $n$. This conjectured holds for graphs of sufficiently large minimum degree~\cite{Bujtas21,BH21}. In the general case, Kostochka and Stodolsky~\cite{KS05} disproved the conjecture by constructing a sequence of connected cubic graphs whose domination ratios $\gamma(G)/n$ converge to a value of at least $1/3+1/69$. The conjecture fails even for $2$-connected cubic graphs. Kelmans~\cite{Kelmans} constructed an infinite family of such graphs with domination ratio $1/3+1/60$, and Stodolsky~\cite{Stodolsky} independently obtained a family with limiting domination ratio at least $1/3+1/78$.

	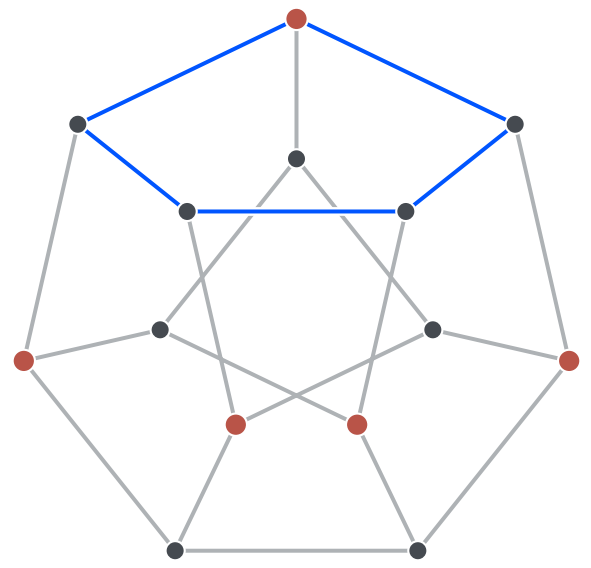
\begin{figure}[htbp]
		\centering
		\definecolor{GPedge}{HTML}{AEB2B5}
		\definecolor{GPvertex}{HTML}{454A50}
		\definecolor{GPcycle}{HTML}{0055FF}
		\definecolor{GPset}{HTML}{B95448}
		
		\begin{tikzpicture}[
			line cap=round,
			line join=round,
			edge/.style={draw=GPedge, line width=0.75pt},
			pentagon/.style={
				edge, draw=GPcycle,
				preaction={draw=white, line width=1.65pt}
			},
			vertex/.style={
				circle, draw=white, fill=GPvertex, line width=0.4pt,
				inner sep=0pt, outer sep=0pt, minimum size=3.6pt
			},
			selected/.style={
				vertex, fill=GPset, minimum size=4.2pt
			}
			]
			\foreach \i in {0,...,6} {
				\coordinate (u\i) at ({90-360*\i/7}:1.85cm);
				\coordinate (v\i) at ({90-360*\i/7}:0.925cm);
			}
			
			\draw[edge] (u0) -- (u1) -- (u2) -- (u3) --
			(u4) -- (u5) -- (u6) -- cycle;
			\draw[edge] (v0) -- (v2) -- (v4) -- (v6) --
			(v1) -- (v3) -- (v5) -- cycle;
			\foreach \i in {0,...,6}
			\draw[edge] (u\i) -- (v\i);
			
			\draw[pentagon] (u6) -- (u0) -- (u1) --
			(v1) -- (v6) -- cycle;
			
			\foreach \p in {u1,u3,u4,u6,v0,v1,v2,v5,v6}
			\node[vertex] at (\p) {};
			\foreach \p in {u0,u2,u5,v3,v4}
			\node[selected] at (\p) {};
		\end{tikzpicture}
		
		\caption{The generalized Petersen graph $P(7,2)$}
		\label{fig:petersen-7-2}
	\end{figure}

	Since these counterexamples contain short cycles, it is natural to ask whether $\gamma(G)\le n/3$ holds for cubic graphs $G$ of order $n$ with larger girth. In this direction, Kawarabayashi, Plummer, and Saito~\cite{KPS} established the
	following bound.
	
	\begin{theorem}[Kawarabayashi, Plummer, and Saito~\cite{KPS}]
		Let $G$ be a connected cubic graph of order $n$ and girth $g$.
		If $G$ has a $2$-factor, then
		\[
		\gamma(G)\le \left(\frac13+\frac{1}{9\lfloor g/3\rfloor+3}\right)n.
		\]
	\end{theorem}
	Kostochka and Stodolsky~\cite{KS09} obtained the bound $\gamma(G)\le (1/3+8/(3g^2))n$ for connected cubic graphs without assuming a $2$-factor. Rautenbach and Reed~\cite{RR} and Kr\'al', \v{S}koda, and Volec~\cite{KSV} further improved the asymptotic upper bounds for cubic graphs of large girth. For cubic graphs with small girth, Fisher, Fraughnaugh, and Seager~\cite{FFS} proved $\gamma(G)\le \lfloor 5n/14\rfloor$ when the girth is at least five. Cho, Choi, Kwon, and Park~\cite{CCKP} proved $i(G)\le 5n/14$ when $G$ is $C_4$-free, where $i(G)$ is the minimum size of an independent dominating set of $G$. In 2008, L\"owenstein and Rautenbach~\cite{LR} proved that every cubic graph $G$ of order $n$ and girth at least $83$ satisfies $\gamma(G)\le n/3$. These results suggest that $1/3$ is a natural benchmark for the domination ratio $\gamma(G)/n$. As reported by Dorbec and Henning~\cite{DH}, Verstra\"{e}te posed the following conjecture.
	
	\begin{conjecture}
		Every cubic graph $G$ of order $n$ and girth at least six satisfies
		$\gamma(G)\le n/3$.
	\end{conjecture}
	
	The generalized Petersen graph $P(7,2)$ (see Fig. 1) gives a counterexample at girth 5, with $14$ vertices and domination number $5>14/3$. Dorbec and Henning~\cite{DH} made substantial progress toward resolving
	Verstra\"{e}te's conjecture by proving the following theorem.
	
	\begin{theorem}[Dorbec and Henning~\cite{DH}]
		Let $G$ be a cubic graph of order $n$ and girth at least six.
		If $G$ contains neither a $7$-cycle nor an $8$-cycle, then $\gamma(G)\le n/3$.
	\end{theorem}
	
	In particular, their theorem gives the $n/3$ bound for every cubic graph of girth at least nine.  In this paper, we prove Verstra\"{e}te's conjecture. 
	
	\begin{theorem}\label{thm:cubic}
		Every cubic graph $G$ of order $n$ and girth at least six satisfies $\gamma(G)\le \left\lfloor n/3\right\rfloor.$
	\end{theorem}
	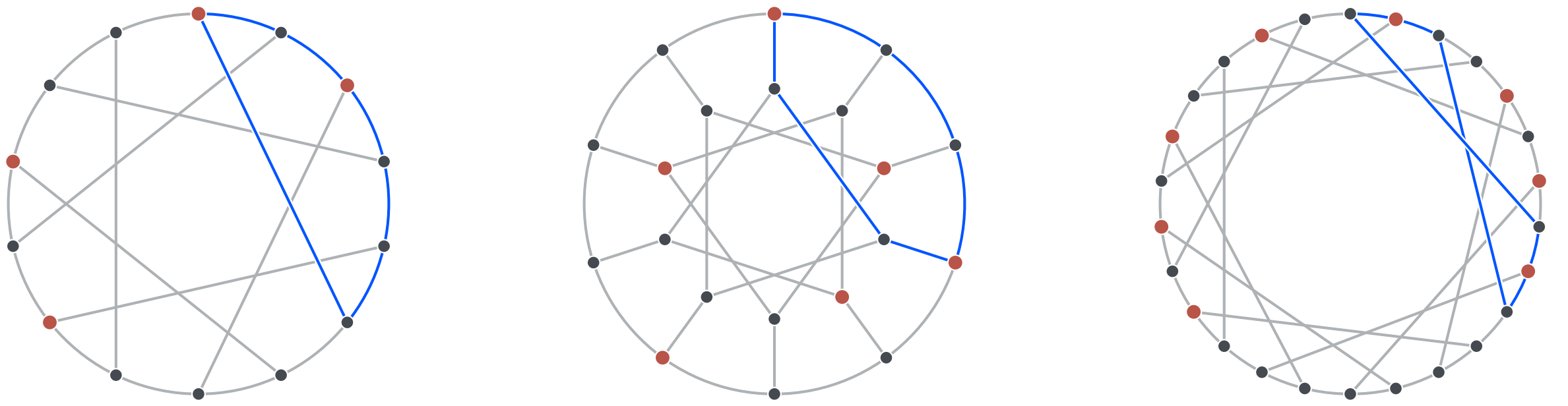
\begin{figure}[htbp]
		\centering
		\definecolor{GPedge}{HTML}{AEB2B5}
		\definecolor{GPvertex}{HTML}{454A50}
		\definecolor{GPcycle}{HTML}{0055FF}
		\definecolor{GPset}{HTML}{B95448}
		
		\resizebox{\linewidth}{!}{%
			\begin{tikzpicture}[
				line cap=round,
				line join=round,
				font=\fontsize{10}{12}\selectfont,
				text=GPvertex,
				edge/.style={draw=GPedge, line width=0.75pt},
				sixcycle/.style={
					edge, draw=GPcycle,
					preaction={draw=white, line width=1.65pt}
				},
				vertex/.style={
					circle, draw=white, fill=GPvertex,
					line width=0.4pt, inner sep=0pt, outer sep=0pt,
					minimum size=3.6pt
				},
				selected/.style={
					vertex, fill=GPset, minimum size=4.2pt
				}
				]
				
				\begin{scope}
					\foreach \i in {0,...,13}
					\coordinate (h\i) at ({90-360*\i/14}:1.85cm);
					
					\draw[edge] (0,0) circle (1.85cm);
					\foreach \i in {0,2,...,12} {
						\pgfmathtruncatemacro{\j}{mod(\i+5,14)}
						\draw[edge] (h\i) -- (h\j);
					}
					
					\draw[sixcycle] (h0)
					arc[start angle=90, end angle={90-360*5/14},
					radius=1.85cm] -- cycle;
					
					\foreach \i in {0,...,13}
					\node[vertex] at (h\i) {};
					\foreach \i in {0,2,9,11}
					\node[selected] at (h\i) {};
					
					\node at (0,-2.56) {(a) Heawood graph};
					\node[font=\fontsize{9}{11}\selectfont] at (0,-3)
					{$n=14,\quad \gamma=i=4$};
				\end{scope}
				
				\begin{scope}[xshift=5.6cm]
					\foreach \i in {0,...,9} {
						\coordinate (u\i) at ({90-36*\i}:1.85cm);
						\coordinate (v\i) at ({90-36*\i}:1.12cm);
					}
					
					\draw[edge] (0,0) circle (1.85cm);
					\foreach \i in {0,...,9} {
						\pgfmathtruncatemacro{\j}{mod(\i+3,10)}
						\draw[edge] (u\i) -- (v\i) (v\i) -- (v\j);
					}
					
					\draw[sixcycle] (u0)
					arc[start angle=90, end angle=-18, radius=1.85cm]
					-- (v3) -- (v0) -- cycle;
					
					\foreach \i in {0,...,9} {
						\node[vertex] at (u\i) {};
						\node[vertex] at (v\i) {};
					}
					\foreach \p in {u0,u3,u6,v2,v4,v8}
					\node[selected] at (\p) {};
					
					\node at (0,-2.56) {(b) Desargues graph $P(10,3)$};
					\node[font=\fontsize{9}{11}\selectfont] at (0,-3)
					{$n=20,\quad \gamma=i=6$};
				\end{scope}
				
				\begin{scope}[xshift=11.2cm]
					\foreach \i in {0,...,25}
					\coordinate (f\i) at ({90-360*\i/26}:1.85cm);
					
					\draw[edge] (0,0) circle (1.85cm);
					\foreach \i in {0,2,...,24} {
						\pgfmathtruncatemacro{\j}{mod(\i+7,26)}
						\draw[edge] (f\i) -- (f\j);
					}
					
					\draw[sixcycle] (f0)
					arc[start angle=90, end angle={90-360*2/26},
					radius=1.85cm]
					-- (f9)
					arc[start angle={90-360*9/26},
					end angle={90-360*7/26}, radius=1.85cm];
					\draw[sixcycle] (f7) -- (f0);
					
					\foreach \i in {0,...,25}
					\node[vertex] at (f\i) {};
					\foreach \i in {1,4,6,8,17,19,21,24}
					\node[selected] at (f\i) {};
					
					\node at (0,-2.56) {(c) $\mathrm{F26A}$};
					\node[font=\fontsize{9}{11}\selectfont] at (0,-3)
					{$n=26,\quad \gamma=i=8$};
				\end{scope}
				
				\node[font=\fontsize{9}{11}\selectfont] at (5.6,-3.67) {
					\textcolor{GPcycle}{\rule[0.5ex]{7mm}{0.75pt}}
					\enspace shortest $6$-cycle
					\qquad
					\textcolor{GPset}{$\bullet$}
					\enspace minimum independent dominating set
				};
				
			\end{tikzpicture}
		}
		
		\caption{Three cubic bipartite graphs of girth six}
	\end{figure}

	Figure~2 shows three graphs that attain the upper bound $\lfloor n/3\rfloor$. In 2009, Kostochka asked whether every cubic bipartite graph $G$ of order $n$ satisfies $\gamma(G)\le n/3$. Kostochka and Stodolsky \cite{KS09} also drew attention to the question, and Henning \cite{Henning} stated it as a conjecture in 2016. A cubic bipartite graph with no $4$-cycle has girth at least $6$. Thus Theorem~\ref{thm:cubic} settles the $C_4$-free case of Kostochka's conjecture.
	
	\begin{corollary}
		Every $C_4$-free cubic bipartite graph $G$ of order $n$ satisfies
		$
		\gamma(G)\le \left\lfloor\frac{n}{3}\right\rfloor.
		$
	\end{corollary}
	
	The remainder of the paper is organized as follows. Section~\ref{sec2} states the weighted bound and proves four estimates for components with leaves. In Section~\ref{sec3}, we narrow down the possible paths of degree-two vertices in a minimum counterexample. In Section~\ref{sec:end}, we rule out the remaining paths and complete the proof by showing that the resulting cubic graph satisfies the bound.

	\section{A weighted bound for marked graphs}\label{sec2}
	
	A graph is \emph{subcubic} if its maximum degree is at most three. For $v\in V(G)$, Let $d_G(v)$ denotes its degree. We use the marked graphs introduced by Dorbec and Henning~\cite{DH}. A mark vindicates that a vertex no longer needs to be dominated after a deletion. A \emph{marked graph} is a subcubic graph in which vertices of degree at most two may carry a mark. A set $D\subseteq V(G)$ is a \emph{marked dominating set} if every unmarked vertex belongs to $D\cup N(D)$. Marked vertices may belong to $D$. Let $\gd(G)$ denote the minimum size of a marked dominating set. Give each marked vertex weight $4$ and each unmarked vertex of degree $0,1,2,3$ weight $12,8,5,4$, respectively. Let $w(G)$ be the total weight of the vertex.
	
	A component is \emph{special} if it has minimum degree at least two and contains a mark. An \emph{exceptional} component is an $8$-cycle with exactly three marks whose successive gaps around the cycle have lengths $2,3,3$. Equivalently, the vertices can be numbered $0,\ldots,7$ in cyclic order with marks at $0,2,5$. Let $c_7(G)$, $s(G)$, and $b_8(G)$ count the unmarked $7$-cycle components, the special components, and the exceptional components, respectively.
	
	\begin{proposition}\label{prop1}
		Every marked subcubic graph $G$ of girth at least six satisfies
		\begin{equation}\label{eq1}
			12\gd(G)\le w(G)+c_7(G)-2s(G)+b_8(G).
		\end{equation}
	\end{proposition}
	
	An unmarked $C_7$ has weight $35$ and domination number three. Its contribution to the right side of~\eqref{eq1} is therefore $35+1=36$. For an exceptional $C_8$, number the vertices as above. The three marks dominate the cycle. Any marked dominating set contains a vertex of $\{0,1,2\}$ to dominate $1$. This vertex dominates at most one of $3,4,6,7$, and any other vertex dominates at most two of them. Hence $\gd(C_8)=3$. Its weight is $3\times 4+5\times 5=37$, so the right side of~\eqref{eq1} is $37-2+1=36$.
	
	If $G$ is cubic, then every vertex is unmarked. Thus $\gd(G)=\gamma(G)$, $w(G)=4|V(G)|$, and $c_7(G)=s(G)=b_8(G)=0$. Proposition~\ref{prop1} gives $12\gamma(G)\le 4|V(G)|$ and therefore implies Theorem~\ref{thm:cubic}.
	
	Suppose Proposition~\ref{prop1} fails and choose a counterexample $G$ of minimum order. Both sides of~\eqref{eq1} are additive over components, so $G$ is connected. An isolated unmarked vertex gives equality, and an isolated marked vertex needs no domination. These graphs and the two cycles considered above satisfy~\eqref{eq1}. Thus $c_7(G)=b_8(G)=0$ and $12\gd(G)>w(G)-2s(G)$, where $s(G)$ is either zero or one.
	
	In each reduction, delete a nonempty set $R\subseteq V(G)$ and retain all marks outside $R$. Let $H$ be the graph obtained from $G-R$ by marking the vertices of a set $M\subseteq V(G-R)$. A cubic vertex can receive a mark only after losing an edge. We use $w(R)$ for the weight of $R$ in $G$ and $w(H)$ for the weight after deletion and marking. For a cycle $Q$ of $G$, write $w(Q)=w(V(Q))$.
	
	Choose a set $S$ of vertices of $G$ that dominates every unmarked vertex of $R$ and every new mark in $H$. Adding $S$ to a minimum marked dominating set of $H$ gives $\gd(G)\le\gd(H)+|S|$. By minimality, $H$ satisfies~\eqref{eq1}. Thus it suffices to show that
	\begin{equation}\label{eq2}
		w(G)-w(H)+2s(H)-b_8(H)-c_7(H)\ge 12|S|+2s(G).
	\end{equation}
	Indeed, these inequalities give
	\[
	12\gd(G)\le 12|S|+12\gd(H) \le 12|S|+w(H)+c_7(H)-2s(H)+b_8(H) \le w(G)-2s(G).
	\]
	
	A special component of $H$ contributes two units to the left side of~\eqref{eq2}, or one if it is exceptional. An unmarked $C_7$ component subtracts one unit.
	
	A vertex of $H$ with a neighbor in $R$ is a \emph{boundary vertex}. On a cycle component of $H$, we call such vertices \emph{attachments}. An attachment of an unmarked $C_7$ was cubic in $G$. If two attachments had a common deleted neighbor, the shorter arc between them and the two edges to that neighbor would form a cycle of length at most five. The same argument applies to two new marks of an exceptional $C_8$, since the shorter arc between any two of its marks has length at most three. In particular, if two new marks share a deleted neighbor and their components have minimum degree at least two, those components contribute at least two units to $2s(H)-b_8(H)$. If the two deleted neighbors are distinct and adjacent, the corresponding cycle may have length six.
	
	A marked leaf $z$ of $H$ that is already dominated by $S$ may also be deleted. If its neighbor $x$ has degree at least two in $H$, then
	\begin{equation}\label{eq3}
		w(H)-w(H-z)=
		\begin{cases}
			4-(8-5)=1, & \text{if $x$ is unmarked and $d_H(x)=2$,}\\
			4-(5-4)=3, & \text{if $d_H(x)=3$,}\\
			4, & \text{if $x$ is marked.}
		\end{cases}
	\end{equation}
	Only the second case can create an unmarked $C_7$, which costs one unit in~\eqref{eq2}. The next estimates apply to smaller components with leaves and allow both old and new marks.
	
	\begin{lemma}\label{lem2}
		Let $F$ be a connected marked subcubic graph of girth at least six with $|V(F)|<|V(G)|$. Suppose every leaf is unmarked, every other vertex has degree two or three, and every marked vertex has degree two.
		\begin{enumerate}[label=\textup{(\roman*)},ref=\roman*]
			\item\label{lem:L1} If $F$ has one leaf, then $12\gd(F)\le w(F)-1$.
			\item\label{lem:L2} If $F$ has one leaf and at most two marks, then $12\gd(F)\le w(F)-2$.
			\item\label{lem:T2} If $F$ is unmarked and has two leaves, each with a cubic neighbor, then $12\gd(F)\le w(F)-2$.
			\item\label{lem:L3} If $F$ has one leaf $U$ and its neighbor $a$ is cubic with an unmarked neighbor $z$ of degree two, then $12\gd(F)\le w(F)-3$.
		\end{enumerate}
	\end{lemma}
	
	\begin{proof}
		Let $U$ be a leaf of $F$ and $a$ its neighbor. In each part we choose $a$ for the dominating set, delete $R=\{U,a\}$, and mark $N(a)\setminus\{U\}$. Call the resulting graph $H$. If $a$ is unmarked of degree two, then $w(R)=8+5=13$. Otherwise $w(R)=8+4=12$. Marking an unmarked neighbor of degree two saves $5-4=1$, and a cubic neighbor keeps weight $4$. Every component of $H$ contains a mark, so $c_7(H)=0$ and $b_8(H)\le s(H)$. By minimality of $G$, we have $12\gd(H)\le w(H)$.
		
		If a vertex $z\in N(a)\setminus\{U\}$ has degree two in $F$, then it is a marked leaf of $H$ dominated by $a$. It may be deleted using~\eqref{eq3} when its remaining neighbor has degree at least two. In part~\textup{(iii)} we verify that this neighbor cannot be the second leaf of $F$. If there are two such vertices $z$, their remaining neighbors are distinct and nonadjacent by girth, so the savings from their deletions add. In each case the choice of $a$ gives $\gd(F)\le 1+\gd(H)$, with $H$ denoting the remainder after any further deletions. Thus the required savings in parts~\textup{(i)}--\textup{(iv)} are $13,14,14,15$, respectively.
		
		\emph{Part \textup{(i)}.} If $w(R)=13$, the first deletion suffices. If $w(R)=12$ and a surviving neighbor of $a$ had degree two, marking it when unmarked or deleting it when already marked gives the extra unit. If every surviving neighbor was cubic, its component now has minimum degree at least two and contains a mark. This also gives at least one unit.
		
		\emph{Part \textup{(ii)}.} If $a$ is unmarked of degree two, then $w(R)=13$. The same argument for the surviving neighbor gives one more unit. Suppose $a$ is marked with $N(a)=\{U,z\}$. If $z$ is cubic, then $H$ is connected, has minimum degree at least two, and has at most two marks. Thus $s(H)=1$ and $b_8(H)=0$, giving $12+2=14$. If $z$ is unmarked of degree two, marking and deleting it saves at least two units.
		
		If $z$ is marked, then $a,z$ are the only marks. Delete $z$ and let $x$ be its remaining neighbor. If $x$ is cubic, this saves at least $3-1=2$ units after the possible $C_7$ correction. If $x$ has degree two in $F$, the deletion saves one unit and leaves $x$ as the only leaf of the remaining connected unmarked graph. Part~\textup{(i)} gives the second unit.
		
		Now let $a$ be cubic, with surviving neighbors $y,z$. If both were cubic, the two new marks have a common deleted neighbor $a$ and give two units. An unmarked neighbor of degree two gives two units when marked and deleted. If both neighbors were marked, deleting the two marked leaves also gives at least two units. The remaining case is that $y$ is cubic and $z$ was marked, up to interchanging them. Delete $z$. If its remaining neighbor is cubic or marked, the saving is at least two units after corrections. Otherwise that neighbor becomes the only leaf of its component. The deletion saves one unit and part~\textup{(i)} gives a second. Any other component contains $y$.
		
		\emph{Part \textup{(iii)}.} If both leaves are adjacent to $a$, delete them with $a$ and mark its third neighbor. The deleted weight is $2\cdot8+4=20$, and $a$ dominates the deleted vertices and the new mark. Otherwise the first deletion saves $12$. If a surviving neighbor $z$ of $a$ had degree two, mark and delete it. Its remaining neighbor cannot be the second leaf, whose neighbor is cubic. Thus the total saving is at least $12+1+1=14$. If both surviving neighbors were cubic, the component containing the second leaf has one leaf and at most two marks. Part~\textup{(ii)} gives the additional two units, and every other component contains a mark.
		
		\emph{Part \textup{(iv)}.} After marking $z$ the saving is $(8+4)+(5-4)=13$. Let $y$ be the third neighbor of $a$. If $y$ is cubic, delete $z$. When the remaining neighbor of $z$ is cubic or marked, the deletion saves at least two more units by~\eqref{eq3}, including the possible $C_7$ correction. Otherwise the deletion saves one unit and leaves that neighbor as the only leaf of its component. Part~\textup{(i)} gives another unit. If $y$ is unmarked of degree two, mark $y$ and delete $z$ for at least $1+1=2$ more units. If $y$ was marked, delete both $y$ and $z$. Their remaining neighbors are distinct, so the two deletions save at least two additional units after corrections. The total saving is at least $15$, giving $12\gd(F)\le w(F)-(15-12)=w(F)-3$.
	\end{proof}
	
	For a disconnected remainder, we can apply Lemma~\ref{lem2} to any component satisfying its hypotheses and~\eqref{eq1} to all other components.
	
	\section{Structure of a minimum counterexample}\label{sec3}
	
	Let $G$ be the minimum counterexample to Proposition~\ref{prop1} fixed above. We use several configurations from Dorbec and Henning~\cite{DH}. When a reduction leaves an unmarked $7$-cycle or an exceptional $8$-cycle, we include its contribution in~\eqref{eq2}.
	
	\begin{claim}\label{claim1}
		The minimum degree of $G$ is at least two. Every marked vertex has degree two, and no two marked vertices are adjacent.
	\end{claim}
	
	\begin{proof}
		Suppose $u$ is an unmarked leaf with neighbor $v$. Take $R=\{u,v\}$, $S=\{v\}$, and $M=N(v)\setminus R$. The vertex $v$ dominates $R\cup M$. Since $G$ has a leaf, $s(G)=0$. Every remaining component contains a mark, so $c_7(H)=0$ and $b_8(H)\le s(H)$. The boundary weights do not increase. Thus $w(G)-w(H)\ge w(R)\ge12=12|S|$, which also covers $G=K_2$.
		
		Now let $u$ be a marked leaf. The preceding case excludes an unmarked leaf as its neighbor. If its neighbor is a marked leaf, then $\gd(G)=0$. Otherwise delete $u$ with $S=M=\emptyset$. If its neighbor is cubic, then $w(G)-w(H)=4-1=3$ and $c_7(H)\le1$. If its neighbor has degree two, then it becomes a leaf, giving $w(G)-w(H)\ge4-3=1$ and $c_7(H)=0$. Both cases give~\eqref{eq2}.
		
		All marks therefore have degree two. Suppose two marks $u,v$ are adjacent, and delete them with $S=M=\emptyset$. Their remaining neighbors are distinct by girth. Each gains at most three in weight, and a cubic neighbor gains only one. Each unmarked $C_7$ in $H$ contains at least one of these cubic neighbors. Thus
		\[
		w(G)-w(H)\ge8-2\cdot3+(3-1)c_7(H)=2+2c_7(H).
		\]
		It follows that $w(G)-w(H)-c_7(H)\ge 2\ge 2s(G)$, as required in~\eqref{eq2}.
	\end{proof}
	
	A cycle is \emph{pendant} if it contains exactly one cubic vertex. Number its vertices in cyclic order starting with $0$ at its cubic vertex, and let $x$ be the neighbor of $0$ off the cycle.
	
	Let $Q$ be a pendant $7$-cycle. If $Q$ has a mark, it has at most three by Claim~\ref{claim1}. Each mark lowers the weight by one, so $w(Q)\ge4+6\cdot5-3=31$. Omit one mark and choose two vertices that dominate the other six consecutive vertices. Use these two for $S$ and delete $Q$ with no new marks. Only $x$ can gain weight, by at most three, and $c_7(H)\le1$. Hence $w(G)-w(H)-c_7(H)\ge31-3-1=27>26\ge12|S|+2s(G)$. If $Q$ is unmarked, take $S=\{x,2,5\}$, delete $Q$ together with $x$, and mark the boundary. Then $M=N(R)$ and $w(G)-w(H)\ge w(R)\ge38$. Both cases give~\eqref{eq2}, so $G$ has no pendant $7$-cycle.
	
	Suppose a pendant $8$-cycle $Q$ has at most two marks. Take $S=\{0,3,5\}$, delete $Q$, and mark $x$. Since $w(Q)\ge4+7\cdot5-2=37$, we need one additional unit. If $x$ was cubic, its component contains a mark and has minimum degree at least two. If $x$ was unmarked of degree two, marking it saves $5-4=1$. If $x$ was already marked, delete it as well. This saves at least one unit. An unmarked $C_7$ arising from this deletion would have been pendant in $G$. In all three cases the corrected decrease is at least $37+1=38\ge12|S|+2s(G)$. Thus every pendant $8$-cycle has at least three marks. The absence of pendant $7$-cycles also implies that an unmarked $C_7$ left by a deletion has at least two attachments.
	
	\begin{claim}\label{cl2}
		No marked vertex has a cubic neighbor.
	\end{claim}
	
	\begin{proof}
		Let $vuw$ be a path with $v$ cubic and $u$ marked. By Claim~\ref{claim1}, $w$ is unmarked. If $w$ is cubic, delete $u$ with $S=M=\emptyset$. Then $w(G)-w(H)=4-1-1=2$. An unmarked $C_7$ in $H$ would have both $v$ and $w$ as attachments. Its shorter $v$--$w$ arc and the path $vuw$ would form a cycle of length at most five. Thus $c_7(H)=0$ and~\eqref{eq2} holds.
		
		Now suppose $w$ has degree two, with other neighbor $x$. Choose $y\in N(x)\setminus\{w\}$ of degree two if possible, and denote the third neighbor of $x$, if any, by $z$. By girth, $v,y,z$ are distinct whenever present, and neither $vy$ nor $vz$ is an edge. If both $x$ and $y$ have degree two, choose $x$ for $S$ and delete $u,w,x,y$ with no new marks. At most one of $x,y$ is marked, so $w(G)-w(H)\ge18-(1+3)=14$. We check the possible $C_7$ correction below.
		
		Suppose every neighbor of $x$ other than $w$ is cubic. Take $S=\{x\}$, delete $u,w,x$, and mark the surviving neighbors of $x$. Each boundary vertex loses one edge, so $\delta(H)\ge2$. Only $v$ can be an unmarked attachment, so $c_7(H)=0$. If $x$ is cubic, the two new marks have the common deleted neighbor $x$. Then $w(G)-w(H)=12$ and $2s(H)-b_8(H)\ge2$. If $x$ is unmarked of degree two, then $w(G)-w(H)=13$ and $2s(H)-b_8(H)\ge1$. In both cases the corrected decrease is at least $14$.
		
		If $x$ is marked, then $w(G)-w(H)=12$. The reduction works unless the component of $y$ is exceptional and is the only special component. This exceptional component also contains $v$, since otherwise it was a pendant $8$-cycle with two marks in $G$. Number its vertices $0,\ldots,7$ with $y=0$ and $v=j$. The path off the cycle is $0-8-9-10-j$, with marks at $8,10$. By girth and symmetry we may take $2\le j\le4$.
		
		The gaps $2,3,3$ and the fact that $j$ is unmarked give the following possibilities, up to reversing the cycle; See Table~\ref{tt1}. In each case the set in the last column dominates all unmarked vertices of $G$.
		Thus $12\gd(G)\le36\le49-2=w(G)-2s(G)$. 
		\begin{table}[h]
			\centering
			\caption{$2\le j\le 4$}\label{tt1}
			\begin{tabular}{@{}lll@{}}
				\toprule
				$j$ & Old marks on the cycle & Dominating vertices \\
				\midrule
				$2$ & $\{3,5\}$ & $\{0,5,10\}$ \\
				$2$ & $\{3,6\}$ & $\{0,4,10\}$ \\
				$3$ & $\{2,5\}$ & $\{0,5,10\}$ \\
				$4$ & $\{2,5\}$ & $\{2,7,10\}$ \\
				$4$ & $\{3,5\}$ & $\{1,6,10\}$ \\
				\bottomrule
			\end{tabular}
		\end{table}

		It remains to consider $d(x)=3$ and $d(y)=2$. Choose $x$ for $S$. If $d(z)=2$, take $R=\{u,w,x,y,z\}$ and $M=\emptyset$. Then $w(G)-w(H)\ge21-(1+3+3)=14$. If $d(z)=3$, write $N(y)=\{x,t\}$ and take $R=\{u,w,x,y\}$ and $M=\{z\}$. When $y$ is unmarked, we have $w(G)-w(H)\ge18-1-3=14$. When $y$ is marked and $d(t)=3$, we have $w(G)-w(H)\ge17-1-1=15$.
		
		If $y$ is marked and $d(t)=2$, apply the cases above to $xyt$ in place of $vuw$. The vertex $r$ beyond $t$ must be cubic with $N(r)=\{t,s,h\}$, where $s$ is marked and $h$ is cubic. Take $R=\{u,w,x,y,t,r,s\}$, $S=\{x,r\}$, and $M=\{z,h\}\setminus R$.
		
		If $|R|=7$, then $w(R)=30$. Only $v$ and the other neighbor of $s$ can gain weight, by at most $1+3$ in total. This bound also holds if they coincide, since a cubic vertex that loses two edges has weight eight. If $rv\in E(G)$, then $h=v$ receives a mark. If $s$ is adjacent to $z$, the mark at $z$ keeps its weight at four even though it loses two edges. All other coincidences are excluded by the degrees or girth. Thus $w(G)-w(H)\ge26$.
		
		If $r=v$, then $s=u$. Take $R=V(vuwxytv)$, a $6$-cycle of weight $26$. Choose $v$ and $x$ for $S$. They dominate the cycle and their outside neighbors. These neighbors are distinct because $G$ has no $5$-cycle. Marking both leaves no unmarked attachment, and neither gains weight.
		
		For the remaining checks, recall that every unmarked $C_7$ needs two attachments. At most three unmarked vertices lie on the boundary, so $c_7(H)\le1$. This is enough if $w(G)-w(H)\ge15$ and $|S|=1$. When the weight bound is $14$ or $26$, an unmarked $C_7$ would use an attachment other than $v$ at which we allowed an increase of three. Such an attachment was cubic and gains only one. The two units recovered cover the correction.
	\end{proof}
	A reduction with no new marks cannot create an exceptional $8$-cycle. Its three old marks would dominate the cycle, and Claim~\ref{cl2} would force every vertex on it to have degree two in $G$. By connectedness, the cycle would already be all of $G$.
	
	Suppose $G=C_n$ and let $m$ be its number of marks. Then $w(G)=5n-m$ and $m\le\lfloor n/2\rfloor$. For $n=6$, two vertices suffice, so $12\gd(G)\le24<25\le w(G)-2s(G)$. On a marked $7$-cycle, omit a mark and choose two vertices to dominate the remaining six vertices. The unmarked $7$-cycle satisfies~\eqref{eq1} by the calculation above. For $n=8$ and at most two marks, choose three vertices. If there are three marks, the gaps are $(2,2,4)$ or $(2,3,3)$. In the first case we may number the marks $0,2,4$ and take $\{2,6\}$. The second case is the exceptional cycle already checked. With four marks, they alternate, so we can number them $0,2,4,6$ and take $\{0,4\}$.
	
	For $n\ge9$, we have $\gd(G)\le\gamma(C_n)=\lceil n/3\rceil$. If $3\mid n$, then $n-m\ge2$ and $12\gd(G)\le4n\le5n-m-2$. If $n\equiv2\pmod3$, then $n\ge11$, so $n-m\ge6$ and $12\gd(G)\le4n+4\le5n-m-2$. Finally suppose $n\equiv1\pmod3$. For $m=0$, we have $12\gamma(C_n)=4n+8\le5n$. For $m>0$, omit a mark and choose $(n-1)/3$ vertices that dominate the remaining path. Then $12\gd(G)\le4n-4\le5n-m-2$. Thus $G$ cannot be a cycle.
	
	A \emph{$2$-path} is a path whose vertices all have degree two in $G$. Since $G$ is connected and is not a cycle, each end of a maximal $2$-path has a cubic neighbor outside the path. These neighbors may coincide.
	
	First, a $2$-path cannot contain six consecutive unmarked vertices. If it does, let $u_1,\ldots,u_6$ be the first six in such a sequence, with $u_0$ the preceding cubic or marked vertex and $u_7$ the vertex after $u_6$. If $u_7=u_0$, the resulting $7$-cycle is either $G$ or a pendant cycle, both of which have been excluded. Choose $u_2,u_5$ for $S$ and delete the six vertices without adding marks. Then $w(G)-w(H)\ge30-1-3=26$. An unmarked $C_7$ in $H$ would require both $u_0,u_7$ as attachments. Both were then cubic, giving $w(G)-w(H)=28$, enough to cover the correction.
	
	Second, there is no path $u_0u_1u_2u_3u_4$ with $u_4$ marked, $u_0$ cubic or marked, and all three internal vertices of degree two. Among the three, only $u_2$ could be marked. Choose $u_2$ and delete $u_1,u_2,u_3,u_4$ without adding marks. Then $w(G)-w(H)\ge18-1-3=14$. By Claim~\ref{cl2}, the neighbor of $u_4$ beyond the path has degree two and becomes a leaf. At most one unmarked attachment remains, so $c_7(H)=0$ and~\eqref{eq2} holds.
	
	\begin{claim}\label{cl3}
		Every maximal $2$-path has at most four vertices.
	\end{claim}
	
	\begin{proof}
		First suppose $u_1,\ldots,u_6$ form a maximal $2$-path, with cubic neighbors $u_0,u_7$ at its ends. The two cubic neighbors are distinct because $G$ has no pendant $7$-cycle. By the preceding observations, only $u_2,u_5$ can be marked. Choose $S=\{u_2,u_5\}$ and delete the six path vertices with no new marks. Then $w(G)-w(H)\ge26$. If an unmarked $C_7$ remains, it contains both $u_0$ and $u_7$, so the path and cycle form all of $G$. If at most one path vertex was marked, then $w(G)-w(H)\ge27$ suffices. If both $u_2,u_5$ were marked, number the cycle so that $u_0,u_7$ are $0,j$, with $j\in\{1,2,3\}$. The set $\{0,j,4\}$ dominates the cycle, and $u_3$ dominates the unmarked vertices on the path. These four vertices give $48\le61-2=w(G)-2s(G)$.
		
		For a maximal $2$-path with at least seven vertices, let $u_0$ be the cubic neighbor at one end, let $u_1,\ldots,u_7$ be its first seven vertices, and let $u_8$ follow $u_7$. If $u_0=u_8$, the two observations above would give a pendant $8$-cycle with at most two marks. Thus $u_0\ne u_8$. At least one of the first seven vertices is marked, and neither $u_1$ nor $u_4$ is marked.
		
		If $u_3$ is marked, then $u_1,u_2,u_4$ are unmarked, and at least two of $u_5,u_6,u_7,u_8$ are unmarked of degree two. Choose $u_1,u_5,u_8$, delete $u_1,\ldots,u_8$, and mark $u_0$ and the surviving neighbors of $u_8$. The deleted weight is at least $37$, and all boundary vertices are marked. If $\delta(H)\ge2$, a special component gives one further unit. Otherwise delete a marked boundary leaf dominated by $S$. An unmarked neighbor of the leaf has degree at least two, since all boundary vertices are marked. Only $u_0$ could have been cubic and lost two edges, which would require $u_0u_8\in E(G)$. In that case its remaining neighbor lost no edge, since $G$ has no triangle. The extra deletion saves at least one unit. It cannot leave an unmarked $C_7$, which would have only one attachment in $G$. Thus the corrected decrease is at least $38=36+2$.
		
		If $u_2$ is marked, take $S=\{u_0,u_4,u_7\}$. When $d(u_8)=2$, delete $u_0,\ldots,u_8$ and mark the surviving neighbors of $u_0$. At least five of the deleted vertices are unmarked of degree two. Only the vertex beyond $u_8$ can gain weight or be an unmarked attachment. Thus $w(G)-w(H)\ge25+16-3=38$ and $c_7(H)=0$. When $u_8$ is cubic, the two observations applied from both ends permit marks only at $u_2,u_6$. Delete $u_0,\ldots,u_7$ and mark $u_8$ and the surviving neighbors of $u_0$. The deleted weight is at least $37$ and all boundary vertices are marked. As in the $u_3$ case, a special component or the deletion of a marked leaf gives the remaining unit.
		
		We may now assume $u_1,\ldots,u_4$ are unmarked. Exactly one of $u_5,u_6$ is marked. Both $u_8$ and the next vertex $u_9$ have degree two. Indeed, if $u_8$ were cubic, the mark would be second or third from that end and the preceding cases would apply. If $u_9$ were cubic, a mark at $u_5$ would be separated from it by three degree-two vertices, contrary to the second observation, and a mark at $u_6$ would be third from that end.
		
		At most two of $u_7,u_8,u_9$ are marked. Take $S=\{u_2,u_5,u_8\}$, $R=\{u_1,\ldots,u_9\}$, and $M=\emptyset$. Then $w(R)\ge42$. The surviving endpoints gain at most four units in total, so $w(G)-w(H)\ge38$. If both are attachments of an unmarked $C_7$, they were cubic and the weight difference is at least $40$. If $u_{10}=u_0$, it becomes a leaf and gains four. These bounds exclude paths on at least seven vertices.
		
		It remains to consider a maximal $2$-path $u_1,\ldots,u_5$ with cubic neighbors $u_0,u_6$ at its ends. Only $u_3$ can be marked. If $u_0=u_6$, the resulting $6$-cycle has weight at least $28$. Choose $u_0,u_3$ and delete the cycle. Let $t$ be the neighbor of $u_0$ off the cycle. If $d(t)=3$, mark $t$. Then $w(G)-w(H)\ge28$ and $c_7(H)=0$. If $d(t)=2$, then $t$ is unmarked by Claim~\ref{cl2}. Delete $t$ as well, giving $w(G)-w(H)\ge28+5-3=30$ and $c_7(H)\le1$.
		
		Suppose $u_0\ne u_6$ and $u_3$ is marked. Delete the five path vertices with $S=\{u_1,u_5\}$ and $M=\{u_0,u_6\}$. Then $\delta(H)\ge2$, $w(G)-w(H)=24$, and $c_7(H)=0$. The marked components contribute at least two units to $2s(H)-b_8(H)$ unless both ends lie on the same exceptional $C_8$. In that case the cycle and the path form $G$. The three marks on the cycle dominate it, and $u_3$ dominates the unmarked vertices on the path that remain to be covered. Thus $12\gd(G)\le48\le61-2=w(G)-2s(G)$.
		
		If all five path vertices are unmarked, take $S=\{u_1,u_4\}$ and mark only $u_0$ after deleting the path vertices. Then $w(G)-w(H)=24$, and $c_7(H)=0$ because only $u_6$ can be an unmarked attachment. The reduction works unless the component of $u_0$ is exceptional and is the only special component. This component also contains $u_6$, since otherwise it was a pendant $8$-cycle with only two marks. The three marks dominate the cycle, and Claim~\ref{cl2} excludes an old mark as a neighbor of the cubic vertex $u_6$. Hence $u_6$ is adjacent to $u_0$. Number the cycle so that $u_0=0$ and $u_6=1$. The old marks are at $3,5$ or at $3,6$. Choosing $0,1$ leaves $4,6$ to be dominated in the first case and $4,5$ in the second. Choose $5$ or $4$, respectively, and add $u_3$ to cover the path. The four vertices again give $48\le61-2=w(G)-2s(G)$.
	\end{proof}
	
	Suppress the degree-two vertices of $G$ to obtain a cubic multigraph, which we call the \emph{compressed graph}. Its vertices are the \emph{roots}. Each edge represents a path in $G$ with no root in its interior. By Claim~\ref{cl3}, these paths have lengths $1,2,3,4$ or $5$. We call the paths and their corresponding edges \emph{short black}, \emph{short green}, \emph{red}, \emph{long black}, and \emph{long green}, respectively.
	
	By Claims~\ref{claim1} and~\ref{cl2}, no mark lies at a root or next to a root, and no two marks are consecutive. Thus paths of length at most three are unmarked. Each longer path has at most one mark. On a long black path it must be the middle vertex.
	
	\begin{claim}\label{cl4}
		If $u b_1 b_2 b_3 v$ is a long black path, then $b_2$ is the only mark in $G$. In particular, there is at most one long black edge.
	\end{claim}
	
	\begin{proof}
		By girth, $u$ and $v$ are distinct. Choose $b_2$ and delete $R=\{b_1,b_2,b_3\}$ with no new marks. Each end gains $5-4=1$. Thus $w(G)-w(H)=14-2=12$ if $b_2$ is marked, and $w(G)-w(H)=15-2=13$ otherwise. We have $\delta(H)\ge2$ and $b_8(H)=0$.
		
		An unmarked $C_7$ in $H$ would contain $u,v$. Together with the deleted path it would form all of $G$, so only $b_2$ could be marked. If $b_2$ were unmarked, then $s(G)=0$ and $w(G)-w(H)-c_7(H)=13-1=12$, as required by~\eqref{eq2}. If $c_7(H)=0$, each component with a surviving mark is special and nonexceptional, and contributes two units. If $G$ has no marks, then $w(G)-w(H)=13$ and $s(G)=0$.
		
		This leaves only the case in which $b_2$ is the sole mark. Two long black paths would each contain this mark in their disjoint interiors, so there is at most one.\end{proof}
	
	\begin{claim}\label{cl5}
		The compressed graph is simple. Each unmarked $C_7$ component left by a deletion has at least three attachments.
	\end{claim}
	
	\begin{proof}
		By Claim~\ref{cl3}, a loop would give a cycle of length at most five in $G$. Suppose two edges of the compressed graph have the same ends. Let the corresponding paths in $G$ have lengths $p,q$ with $p\le q$. Then $p+q\ge6$, $q\le5$, and $(p,q)\ne(4,4)$ by Claim~\ref{cl4}. Their union is a cycle $Q$. Number its vertices in cyclic order from $0$ to $p+q-1$, with the two roots at $0,p$.
		
		First take $R=V(Q)$ and $M=N(R)$. For $(p,q)=(2,4)$, the mark is at $4$, so $\{0,2\}$ suffices and $w(Q)=27$. For $(3,3)$ and $(3,4)$, take $\{0,3\}$, with $w(Q)=28$ and $32$, respectively. In the latter case the mark is at $5$. For $(4,5)$, take $\{0,4,6\}$, since $2$ is the only mark and $w(Q)=42$. Each choice includes the two roots and therefore dominates $M$. The bound $w(Q)\ge12|S|+2s(G)$ gives~\eqref{eq2}.
		
		For $(p,q)=(3,5)$, take $S=\{0,3,5\}$. At most one of $5,6$ is marked, so $w(Q)\ge37$. The roots have distinct outside neighbors, since a common neighbor and the path of length three would give a $5$-cycle. If both outside neighbors were cubic, their marked components have minimum degree at least two, giving $2s(H)-b_8(H)\ge1$. Otherwise an outside neighbor had degree two and was unmarked by Claim~\ref{cl2}. Marking it saves $5-4=1$. The corrected decrease is therefore at least $37+1=38$.
		
		For $(p,q)=(1,5)$, write the cycle as 
		$u v y_1 y_2 y_3 y_4 u$ and take $S=\{u,y_2\}$. We have $w(Q)\ge27$. Let $a,b$ be the outside neighbors of $u,v$, respectively. If $d(a)=2$, then $a$ is unmarked. Delete $a$ with $Q$ and add no marks. Its remaining neighbor $a'$ differs from $b$, since $a'=b$ would give a $4$-cycle. Then $w(G)-w(H)\ge32-3-3=26$. An unmarked $C_7$ would use both $a'$ and $b$ as attachments, so both were cubic and $w(G)-w(H)\ge30$. Reversing the labeling gives the same reduction when $d(b)=2$. If $d(a)=d(b)=3$, take $R=V(Q)$ and $M=\{a\}$. Here $w(G)-w(H)\ge27-1=26$ and only $b$ can be an unmarked attachment.
		
		For $(p,q)=(2,5)$, the cycle has order seven. A mark can lie only at $4$ or $5$. If it is at $4$, take $S=\{2,6\}$. If it is at $5$, take $S=\{0,3\}$. In each case take $R=V(Q)$ and mark the outside neighbor of the root in $S$. Then $w(G)-w(H)\ge32-3=29$ and $c_7(H)\le1$.
		
		If this $7$-cycle is unmarked, let $a$ be the outside neighbor of $0$. Take $R=V(Q)\cup\{a\}$, $M=N(R)$, and $S=\{a,2,5\}$. Since $w(Q)=33$, the deleted weight is $33+5=38$ if $d(a)=2$ and $33+4=37$ if $d(a)=3$. A degree-two vertex in $M$ was unmarked, so marking it saves $5-4=1$. If all boundary vertices were cubic, two edges leaving $R$ cannot have the same end in $H$, since that would give a cycle of length at most five. Then $\delta(H)\ge2$ and $2s(H)-b_8(H)\ge1$.
		
		Finally, for $(p,q)=(5,5)$, write $Q=u x_1 x_2 x_3 x_4 v y_1 y_2 y_3 y_4 u$ and choose $u,x_3,y_2$. Delete all vertices of $Q$ except $v$ and mark the neighbor of $u$ off $Q$. There are at most two marks on $Q$, so $w(R)\ge42$. By the $(1,5)$ case, $uv\notin E(G)$, so $d_H(v)=1$ and its weight rises from four to eight. Only $v$ can be an unmarked attachment. Thus $w(G)-w(H)\ge38$ and $c_7(H)=0$. This calculation also permits a third path between $u$ and $v$.
		
		Now consider the attachments of an unmarked $C_7$ in $H$. A single attachment would make the cycle pendant in $G$. With two attachments, the arcs between them have lengths $(1,6)$, $(2,5)$, or $(3,4)$. In the first case the longer arc has five internal degree-two vertices, contrary to Claim~\ref{cl3}. The second case gives parallel edges in the compressed graph. In the third, the longer arc is an unmarked long black path, contrary to Claim~\ref{cl4}. Therefore there are at least three attachments.
	\end{proof}
	
	By~\eqref{eq3}, deleting a marked leaf $z$ already dominated by $S$ saves at least one unit after correction when its neighbor has degree at least two in $H$. If $z$ is the only leaf and all other vertices have degree two or three, the saving is at least two. Indeed, when the deletion saves only one unit, the neighbor of $z$ becomes the only leaf of its component and \leafpart{lem:L1} gives another unit. An isolated mark has weight $4$ and can be deleted for that saving. If two marked leaves have a common deleted neighbor, girth makes their surviving neighbors distinct and nonadjacent. Under the same degree assumption, their deletions save at least $1+1=2$.
	
	\begin{claim}\label{cl6}
		No long black edge is adjacent to a green or red edge. The green edges form a matching. There is no green--red--green path, and no red edge is adjacent to a long green edge.
	\end{claim}
	
	\begin{proof}
		Let a green path $P$ join $v$ to $w$ at one end of the long black path $u b_1 b_2 b_3 v$. By Claim~\ref{cl4}, $b_2$ is the only mark. Delete $R=\{u,b_1,b_2,b_3\}\cup(V(P)\setminus\{w\})$ and mark the surviving neighbors of $u,v$. If $P$ is short, take $S=\{u,v\}$, giving $w(R)=27$. If $P=v a_1 a_2 a_3 a_4 w$ is long, take $S=\{u,v,a_3\}$, giving $w(R)=42$. Only $w$ can gain weight, and it gains at most one because it receives a mark if $uw\in E(G)$. Thus $c_7(H)=0$ and
		\[
		w(G)-w(H)\ge w(R)-1\ge12|S|+2.
		\]
		
		Suppose instead that $v r_1 r_2 w$ is red. Take 
		$R=\{u,v,w,b_1,b_2,b_3,r_1,r_2\},$ $M=N(R),$ and $S=\{u,v,w\}.$ 
		The set $S$ dominates every unmarked vertex of $R$ and every new mark. Since $w(R)=36=12|S|$ and all boundary vertices are marked, we need two further units. If a boundary vertex had degree two in $G$, marking it saves $5-4=1$ and leaves it with degree at most one. We can delete it for at least one more unit. If it has an unmarked neighbor, that neighbor lies outside the boundary and still has degree at least two, so~\eqref{eq3} applies.
		
		Now suppose all boundary vertices were cubic. A leaf of $H$ must have been adjacent to $u,v$ or to $u,w$, since a common neighbor of $v,w$ would give a $5$-cycle. Each pair has at most one common neighbor. If $H$ has two leaves, they share the deleted neighbor $u$ and can be deleted for two units. If it has exactly one leaf and no isolated vertex, its deletion saves two units by the single-leaf estimate above. Deleting an isolated mark saves four units. We may therefore assume $\delta(H)\ge2$.
		
		If $uw\notin E(G)$, the marks at the two surviving neighbors of $u$ give two units. If $uw\in E(G)$, let $a,b,c$ be the remaining neighbors of $u,w,v$, respectively. The three marks are distinct and give at least two units unless all three lie on one exceptional $8$-cycle. Since each remaining component contains a mark, this cycle is then all of $H$. The distance between $a$ and $b$ on it is three. A shorter distance would give a cycle of length at most five with $a u w b$, and the gaps $2,3,3$ give the upper bound. Number the cycle with $a=0$ and $b=3$, so $c=5$ or $c=6$. Write $(u,v,w)=(8,9,10)$. The two other paths are $8-11-12-13-9$ and $9-14-15-10$, with a mark at $12$. For $c=5$ take $\{0,2,5,12,15\}$, and for $c=6$ take $\{0,3,5,12,14\}$. Both choices dominate all unmarked vertices in the respective graphs, so $12\gd(G)\le60<73-2=w(G)-2s(G)$.
		
		Now suppose two green edges share a root $u$. Let $v,w$ be their other ends and $t$ the neighbor of $u$ off their paths. By Claim~\ref{cl5}, $v\ne w$. For $S$, choose $u$ and the third internal vertex from $u$ on each long green path. These vertices dominate $t$ and both interiors. A short green path has internal weight $5$. A long green path has internal weight at least $4\times 5-1=19$, since it has at most one mark. Thus each of the $|S|-1$ long paths contributes at least $14$ more than a short path.
		
		If $d(t)=2$, then $t$ is unmarked. Delete $t,u$ and the two interiors with $M=\emptyset$. The other neighbor $z$ of $t$ is distinct from $v,w$ by simplicity. Hence
		\[
		w(G)-w(H)\ge19+14(|S|-1)-(1+1+3)=14|S|.
		\]
		If $z$ is marked or becomes a leaf, only $v,w$ can be unmarked attachments, so $c_7(H)=0$. If $z$ was cubic, it gains one instead of three and $c_7(H)\le1$. In either case $w(G)-w(H)-c_7(H)\ge14|S|\ge12|S|+2$.
		
		If $d(t)=3$, retain $t$ and mark it. Then $d_H(t)=d_H(v)=d_H(w)=2$ and $\delta(H)\ge2$. Only $v,w$ can be unmarked attachments, so $c_7(H)=0$. The same count gives $w(G)-w(H)\ge14+14(|S|-1)-2=14|S|-2$. If $|S|\ge2$, this is at least $12|S|+2$. If $|S|=1$, both green paths are short. The component of $t$ cannot be exceptional. On such a cycle, any root $v$ or $w$ would be adjacent to a mark because the three marks dominate the cycle. By Claim~\ref{cl2}, this mark would have to be $t$, giving a $4$-cycle through $u$. If neither root lies on the cycle, it was pendant in $G$ with only two marks. Hence the component of $t$ gives the remaining two units. This proves that the green edges form a matching.
		
		Suppose $u_1u_2u_3u_4$ is a green--red--green path in the compressed graph. Simplicity and the green matching ensure that the four roots are distinct. Delete $u_2,u_3$ and the interiors of the three paths. Mark the surviving third neighbors of the deleted roots. For $S$, choose $u_2,u_3$ and the third internal vertex of each long green path, counted from its deleted root. Only $u_1,u_4$ can gain weight. An end adjacent to the opposite deleted root receives a mark and retains weight $4$. Otherwise it loses one edge and gains one unit. With two short green paths the deleted weight is $28$. Each of the $|S|-2$ long paths adds at least $14$, so
		\[
		w(G)-w(H)\ge28+14(|S|-2)-2=14|S|-2\ge12|S|+2.
		\]
		There are at most two unmarked attachments, so $c_7(H)=0$.
		
		Finally, suppose $u r_1 r_2 v$ is red and $v a_1 a_2 a_3 a_4 w$ is long green. Take
		$R =\{u,v,r_1,r_2,a_1,a_2,\allowbreak  a_3, a_4\},$
		$M =(N(u)\cup N(v))\setminus R,$ and  $S=\{u,v,a_3\}.$ 
		Only $w$ can gain weight, so $w(G)-w(H)\ge w(R)-1\ge37-1=36$. By the preceding restrictions, any vertex of $M$ that had degree two is the first internal vertex of a red path. Marking and deleting it saves at least $1+1=2$ units, and deleting it saves more if it is isolated. The second internal vertex becomes a leaf, so the extra deletion cannot create an unmarked $C_7$.
		
		Suppose that all vertices in $M$ were cubic. If $uw\notin E(G)$, every boundary vertex loses exactly one edge. Indeed, a common neighbor of $u,v$ would produce a $5$-cycle, and the only other possible vertex that could lose two edges is $w$. The two marks at the surviving neighbors of $u$ then give two units. If $uw\in E(G)$, then $w\in M$ and $d_H(w)=1$, so $w(G)-w(H)\ge37$. Its remaining neighbor, if unmarked, has degree at least two because every boundary vertex is marked. Delete $w$ for another unit. At most two unmarked attachments remain, so the corrected decrease is at least $38=12|S|+2$ in both cases.
	\end{proof}
	
	\begin{claim}\label{cl7}
		There is no marked long green path, and the compressed graph contains no all-black triangle.
	\end{claim}
	
	\begin{proof}
		Suppose the long green path $P=u a_1 a_2 a_3 a_4 v$ has a mark. By reversing $P$ if necessary, we may assume $a_3$ is marked. By Claim~\ref{cl6}, all other edges incident with $u$ or $v$ are short black. Take $R=V(P)$, $M=N(v)\setminus R$, and $S=\{a_1,v\}$. Then $w(R)=27$. Only the two neighbors of $u$ off $P$ can be unmarked attachments, so $c_7(H)=0$. If $N(u)\cap N(v)=\emptyset$, each gains one unit, giving $w(G)-w(H)\ge27-2=25$. A special component gives one more. If $u,v$ have a common neighbor, it is unique by girth and belongs to $M$, so $w(G)-w(H)\ge27-1=26$.
		
		Suppose the compressed graph has an all-black triangle. Exactly one of its edges is long, since three short edges would give a triangle in $G$ and Claim~\ref{cl4} permits at most one long black edge. Let $u b_1 b_2 b_3 v$ be its long path, and let $uw,vw$ be the other edges. The mark $b_2$ is the only one in $G$. Denote the remaining neighbors of $u,v,w$ by $x,e,t$, respectively. Both $x$ and $e$ are cubic. If $d(t)=3$, take $R=\{u,v,w,b_1,b_2,b_3\}$, $M=\{x,e\}$, and $S=\{u,v\}$. The vertices $x,e,t$ are distinct because a coincidence would give a cycle of length three or four. Then $w(G)-w(H)=26-1=25$, and a special component gives at least one more unit. Only $t$ is a possible unmarked attachment.
		
		If $d(t)=2$, include $t$ in $R$, take $S=\{w,b_2\}$, and put $M=\emptyset$. The neighbor $z$ of $t$ other than $w$ is unmarked. By girth, $x,e,z$ are distinct. Now $w(R)=31$. If $d(z)=2$, then $d_H(z)=1$, so $w(G)-w(H)=31-1-1-3=26$ and only $x,e$ can be unmarked attachments. If $d(z)=3$, then $w(G)-w(H)=28$ and $c_7(H)\le1$. Both cases give~\eqref{eq2}.
	\end{proof}
	
	\begin{claim}\label{cl8}
		There is no long black path.
	\end{claim}
	
	\begin{proof}
		Suppose $u b_1 b_2 b_3 v$ is a long black path. Its middle vertex $b_2$ is the only mark in $G$. All other edges at $u$ or $v$ are short black. Let $N(u)\setminus\{b_1\}=\{w,x\}$ and $N(v)\setminus\{b_3\}=\{d,e\}$. A common vertex would give an all-black triangle, contrary to Claim~\ref{cl7}. The following restrictions on the edges at $w$ also hold at $x$ by symmetry.
		
		Suppose all three edges at $w$ are black. Take $R=\{u,v,w,b_1,b_2,b_3\}$, $S=\{u,v\}$, and $M=(N(u)\cup N(v))\setminus R$. Each neighbor of $w$ other than $u$ can gain one unit, giving $w(G)-w(H)\ge26-2=24$. If $\delta(H)\ge2$, the marks $d,e$ give two further units. Otherwise a common neighbor $y$ of $v,w$ is a marked leaf. It is unique by girth, and any other coincidence would give a cycle of length at most five or an all-black triangle. Since $y$ is marked, it gains no weight, so $w(G)-w(H)\ge25$. An unmarked boundary vertex was cubic and lost only one edge. Hence an unmarked neighbor of $y$ has degree at least two in $H$. Deleting $y$ saves at least one more unit by~\eqref{eq3} and leaves at most two unmarked attachments.
		
		Now let $ww'$ be a green edge in the compressed graph and $t$ the third neighbor of $w$ in $G$. If the path is $w g w'$, take
		$ R=\{u,w,b_1,b_2,b_3,g\},~ M=\{t\},$ and $ S=\{w,b_2\}.$
		The vertices $v,x,w'$ are distinct. Indeed, $w'\ne v$ because a green edge cannot meet the long black edge, and $w'\ne x$ by the absence of a $4$-cycle. Also, $v\ne x$ by the exclusion of all-black triangles. These three vertices are outside $M$ and each loses one edge. We have $w(R)=27$. If $d(t)=3$, then $t\notin\{v,x,w'\}$ by simplicity and the same exclusions. Hence $w(G)-w(H)=24$ and the only mark belongs to a special nonexceptional component. If $d(t)=2$, then $wt$ starts a red path. Marking and deleting $t$ gives two additional units, and the second internal vertex becomes a leaf.
		
		An unmarked $C_7$ would have attachments $v,x,w'$. Up to symmetry the arc lengths between them are $(1,1,5)$, $(1,2,4)$, $(1,3,3)$, or $(2,2,3)$. The pair of length-two arcs in $(2,2,3)$ would give adjacent green edges. The length-four arc in $(1,2,4)$ would give an unmarked long black path, and one of the length-three arcs in $(1,3,3)$ would give a red edge at $v$. In the case $(1,1,5)$, the length-one arcs must meet at $v$ because no green edge meets the long black edge. Then $vx$ is an edge, and $u,v,x$ form an all-black triangle. Hence $c_7(H)=0$.
		
		Now suppose $w a_1 a_2 a_3 a_4 w'$ is long green. This path is unmarked and $d(t)=3$. Take $R=\{u,w,b_1,b_2,b_3,a_1,a_2,a_3,a_4\}$, $M=\{t\}$, and $S=\{w,b_2,a_3\}$. Suppose $w'\ne x$. Simplicity and the long black edge at $v$ give $w'\notin\{u,v\}$ and $t\ne w'$. The equality $t=v$ would give an all-black triangle, and $t=x$ would give a triangle in $G$. Thus $w(R)=42$, and each of $v,x,w'$ gains one unit. We obtain $w(G)-w(H)=42-3=39$. Since $c_7(H)\le1$, the corrected decrease is at least $38$.
		
		If $w'=x$, let $Q=u w a_1 a_2 a_3 a_4 w' u$. Take $R=V(Q)\cup\{b_1,b_2,b_3,v\}$, $M=\{d,e\}$, and $S=\{u,v,a_2,a_4\}$. We have $w(R)=50$. The outside neighbors $f,h$ of $w,w'$ are cubic and distinct, since a common neighbor would give a $4$-cycle through $u$. An edge from $v$ to either $w$ or $w'$ would give an all-black triangle. Only $f,h$ can gain weight. Each gains one unit unless it belongs to $M$, in which case it becomes a marked leaf. Thus $w(G)-w(H)=48+|\{f,h\}\cap M|$. If neither lies in $M$, every boundary vertex loses one edge. The marks $d,e$ then give two further units. If both lie in $M$, the difference is already $50$. If exactly one lies in $M$, delete that marked leaf. Its deletion saves at least one unit, since an unmarked boundary neighbor was cubic and lost only one edge. At most two unmarked attachments remain.
		
		Suppose two red paths $w r s A$ and $w p q B$ meet at $w$. Take $R=\{u,w,b_1,b_2,b_3,r,p\}$, $S=\{w,b_2\}$, and $M=\emptyset$. The vertices $s,q$ become leaves, and $v,x$ each lose one edge. Thus $w(G)-w(H)=32-2\times 3-2=24$. The neighbors $A,B$ of the leaves are still cubic. Indeed, $A=B$ would give parallel edges, and an end in $\{u,w,x\}$ would give a short cycle or parallel edges. Neither end can be $v$, where the long black edge excludes a red edge.
		
		If $s,q$ belong to the same component, apply \leafpart{lem:T2} to it. If they lie in different components, apply \leafpart{lem:L1} to each. Only $v,x$ can be attachments of an unmarked $C_7$ in another component, so there is no such component. In either case $12\gd(H)\le w(H)-2$, giving
		\[
		12\gd(G)\le12|S|+w(H)-2=24+w(H)-2=w(G)-2.
		\]
		
		Thus $w,x$ each have one red path, say $w r_1 r_2 w'$ and $x s_1 s_2 x'$. Let $f,h$ be the short black neighbors of $w,x$ other than $u$, respectively. Take
		$R =\{u,v,w,x,b_1,b_2,b_3,r_1,r_2,s_1,s_2\},$
		$M =\{w',x',d,e\},$ and $S=\{u,v,r_2,s_2\}.$ 
		The set $S$ dominates the unmarked vertices of $R$ and the new marks, and $w(R)=50$. By girth, $f,h$ are distinct cubic vertices. They are the only vertices that can gain weight. As above,
		$w(G)-w(H)=48+|\{f,h\}\cap M| $ and $ c_7(H)=0.$
		If $\delta(H)\ge2$, the marks $d,e$ share the deleted neighbor $v$ and give two more units.
		
		To determine which boundary vertices can become leaves or isolated vertices, first consider $w',x'$. Each loses the edge on its red path and must lose at least one more. By girth it cannot be adjacent to $u$ or to the starting root of its red path. An edge to $v$ is possible. An edge to the other starting root gives $w'=h$ or $x'=f$ and a $6$-cycle through $u$. The two far ends may also coincide. A common neighbor of $v,w$ or of $v,x$ similarly belongs to $\{f,h\}\cap M$. These are all possibilities, since every deleted internal vertex has both neighbors on its path.
		
		If $f,h\in M$, then $w(G)-w(H)=50$. If exactly one lies in $M$, the difference is $49$. Delete this marked leaf or isolated vertex, which is dominated by $v$ or by the chosen internal vertex on its red path. Its deletion saves at least one unit, since an unmarked boundary neighbor was cubic and lost only one edge. If $w'=x'$ is a second marked leaf with the same degree-two neighbor, the calculation still applies. We delete only the vertex in $\{f,h\}\cap M$, leaving the common neighbor as a leaf. At most two unmarked attachments remain.
		
		Suppose $\{f,h\}\cap M=\emptyset$ and $w'=x'$. If this common end meets $v$, it is isolated in $H$, so its deletion saves $4$. Otherwise it is the only leaf. Its remaining neighbor differs from $f,h$, since equality would give a $5$-cycle. The single-leaf estimate then gives two units on deleting it.
		
		It remains that $w'\ne x'$ and at least one, say $w'$, is adjacent to $v$. Change the reduction to
		$ R =\{u,v,w,w',b_1,b_2,b_3,r_1,r_2\},$
		$M =(N(w)\cup N(w'))\setminus R$ and  $ S=\{w,w',b_2\}.$
		Now $w(R)=40$. Only $x$ and the remaining neighbor of $v$ can gain weight, by at most one unit each. If that neighbor also meets $w$, it receives a mark. The same holds for $x$ if it meets $w'$. All other coincidences would give an all-black triangle or a cycle of length at most five. Thus $w(G)-w(H)\ge38=12|S|+2$. At most two unmarked attachments remain, so~\eqref{eq2} holds.
	\end{proof}
	
	Only long black and long green paths could carry marks. By Claims~\ref{cl7} and~\ref{cl8}, $G$ is therefore unmarked and $s(G)=0$.
	
	Suppose $u a_1 a_2 a_3 a_4 v$ is long green. All other edges at $u,v$ are short black. Take $R=\{a_1,a_2,a_3,a_4,v\}$, $M=\{u\}$, and $S=\{a_1,a_4\}$. Simplicity excludes $uv$, so every boundary vertex loses one edge. The two surviving neighbors of $v$ each gain one unit and $u$ retains weight $4$. We have $\delta(H)\ge2$, so $s(H)=1$ and $b_8(H)=0$. Since only two unmarked attachments are available, $c_7(H)=0$. Thus $w(G)-w(H)+2s(H)=(24-2)+2=24=12|S|$, as required in~\eqref{eq2}. The remaining edges are short black, short green, or red.
	
	\section{Completion of the proof}\label{sec:end}
	
	\begin{proof}[Proof of Proposition~\ref{prop1}]
		By Section~\ref{sec3}, $G$ is connected and unmarked, and its compressed graph is simple. Call paths of lengths one, two, and three black, green, and red, respectively. The green edges form a matching, and there is no green--red--green path. Call a cubic vertex on a cycle of $G$ a \emph{port} of the cycle. For an unmarked $C_7$ component of $H$, the ports are exactly the attachments. By Claim~\ref{cl5}, there are at least three. Each has one neighbor in $R$, and these neighbors are distinct by girth. Since $s(G)=0$, the right side of~\eqref{eq2} is $12|S|$.
		
		\smallskip
		\noindent\refstepcounter{finishstep}\label{lem:no-redgreen}\emph{Step \thefinishstep. No red edge is adjacent to a green edge.}\par\nopagebreak[4]
		
		Suppose the green path $U h V$ and the red path $V r s W$ meet at $V$. If $UW\in E(G)$, the paths and $UW$ form a $6$-cycle $Q$. The outside neighbors of $U,V,W$ are distinct by girth. Take $S=\{V,W\}$, delete $R=V(Q)$, and mark the outside neighbors of $V,W$. The set $S$ dominates $R\cup M$. Only the outside neighbor of $U$ can gain weight or be an unmarked attachment. Hence $c_7(H)=0$ and $w(G)-w(H)\ge27-3=24=12|S|$.
		
		If $UW\notin E(G)$, take $R=\{h,V,r,s,W\}$ and $S=\{V,W\}$. Mark the surviving neighbor of $V$ and the two surviving neighbors of $W$. These three vertices and $U$ are distinct because a coincidence would give a short cycle or parallel edges in the compressed graph. The deleted weight is $23$. Only $U$ gains a unit, and marking a vertex of degree two in $G$ saves one unit. Thus
		\[
		w(G)-w(H)=22+\bigl|\{x\in M:d(x)=2\}\bigr|.
		\]
		If all vertices of $M$ were cubic, then $\delta(H)\ge2$. The two marks with common deleted neighbor $W$ give $2s(H)-b_8(H)\ge2$. Since only $U$ can be an unmarked attachment, $c_7(H)=0$.
		
		Otherwise, delete a newly marked leaf already dominated by $S$. It starts a red path, since a green path would give adjacent green edges or a green--red--green path. Its second internal vertex still has degree two, since a coincidence with another affected internal vertex would give a loop or parallel red edges in the compressed graph. The first deletion saves at least $22+1=23$, and deleting the leaf saves $4-3=1$. The second vertex becomes a leaf, so $U$ is still the only possible unmarked attachment. The total saving is at least $23+1=24$.
		
		\smallskip
		\noindent\refstepcounter{finishstep}\label{lem:red-matching}\emph{Step \thefinishstep. The red edges form a matching, and each $7$-cycle has at least four ports.}\par\nopagebreak[4]
		
		Let $S$ be the root set of a connected subgraph of the compressed graph consisting of $e$ red edges. In $G$, the roots in $S$ dominate the $2e$ internal vertices of those red paths and all neighbors outside the subgraph. Delete the $|S|+2e$ vertices and mark the boundary. Then $c_7(H)=0$ and $w(G)-w(H)\ge4|S|+10e$, so $5e\ge4|S|$ suffices for~\eqref{eq2}. This holds on a red cycle, where $e=|S|$, and on a red tree with five roots, where $e=4$. Every larger red tree contains one on five roots. Thus the red edges form a forest whose components have at most four roots.
		
		Suppose a component is a star on four roots, with paths $v a_i b_i v_i$ for $1\le i\le3$. By Step~\ref{lem:no-redgreen} and the bound of four roots, each $v_i$ has two black edges besides its red edge. For outer roots labeled $v_1,v_2$, use
		$R =\{v,v_1,v_2,a_1,b_1,a_2,b_2,a_3,b_3\},$ 
		$ M =\emptyset, $ and $S=\{b_1,b_2,a_3\}.$ 
		The set $S$ dominates $R$, and $w(R)=42$. If two outer roots are adjacent, label them $v_1,v_2$. Three edges then join $R$ to $H$, all with cubic ends in $G$. Only $v_3$ can be the end of two such edges. One of these edges is $b_3v_3$, and $v_3$ meets at most one of $v_1,v_2$ because $G$ has no triangle. Allowing one unit per edge and two more if $v_3$ becomes a leaf gives $w(G)-w(H)\ge42-3-2=37$.
		
		Suppose the outer roots are pairwise nonadjacent. If a pair has no common neighbor, label it $v_1,v_2$. The five edges leaving $R$ have distinct cubic ends, so $w(G)-w(H)=42-5=37$. Both reductions leave at most five edges from $R$ to $H$, and each unmarked $C_7$ needs at least three of them. Thus $3c_7(H)\le5$, and $w(G)-w(H)-c_7(H)\ge37-1=36$.
		
		It remains that the outer roots are pairwise nonadjacent and every pair has a common neighbor, unique by girth. If one vertex $z$ meets all three, choose $v,v_1,v_2,v_3$. Delete $z$ with the ten vertices of the star in $G$ and mark the boundary. The chosen roots dominate the deleted vertices and the boundary, and the weight decrease is at least $50$. Otherwise let $z_{ij}$ be the common neighbor of $v_i,v_j$. The outer roots and $z_{12},z_{23},z_{31}$ form a $6$-cycle. Choose the same four roots, delete the cycle together with $v$ and the six internal vertices of the red paths, and add no marks. Then $w(R)=7\cdot4+6\cdot5=58$. The third edge at each $z_{ij}$ leaves $R$ by girth and uniqueness. These three edges have distinct ends in $H$, since each pair of the $z_{ij}$ already shares an outer root. Thus $w(G)-w(H)\ge58-9=49$ and $c_7(H)\le1$. This excludes the star.
		
		Consider a $7$-cycle with arc lengths $(3,3,1)$ between its ports. Let $A r s B$ and $B t u C$ be its red paths, with $CA$ black. The star case and Step~\ref{lem:no-redgreen} imply that the third edge at $B$ is black. Take $R=\{r,s,B,t,u\}$, $M=\{A\}$, and $S=\{r,t\}$. Each unmarked boundary vertex gains one unit. The component of $A$ has minimum degree at least two and a single mark. It also contains $C$, so only the third neighbor of $B$ can attach to an unmarked $C_7$. Hence $s(H)=1$, $b_8(H)=c_7(H)=0$, and $w(G)-w(H)+2s(H)=(24-2)+2=24$. This excludes $(3,3,1)$. The only other arc lengths with three ports are $(3,2,2)$, which would give adjacent green edges. Therefore every $7$-cycle has at least four ports.
		
		If a red component still has at least two edges, it is a path. Let $u_1u_3,u_3u_4$ be consecutive edges with $u_4$ an end of the path. The other edges at $u_3,u_4$ are black. Denote them by $u_2u_3,u_4u_5,u_4u_6$. If $u_1$ equals $u_5$ or $u_6$, the two red paths and one black edge give the excluded arc lengths $(3,3,1)$. Any other coincidence would give a short cycle or parallel edges. Thus the six roots are distinct.
		
		Delete $u_3,u_4$ and the four internal vertices of the two red paths with no new marks. Choose the internal vertex adjacent to $u_3$ on the first path and the one adjacent to $u_4$ on the second. They dominate $R$. We have $w(R)=28$, and each of the four cubic boundary vertices $u_1,u_2,u_5,u_6$ loses one edge. Thus $w(G)-w(H)=28-4=24$. An unmarked $C_7$ would need all four vertices as ports, but $u_5,u_6$ share the deleted neighbor $u_4$. The resulting reduction proves that the red edges form a matching.
		
		\smallskip
		\noindent\refstepcounter{finishstep}\label{lem:no-green}\emph{Step \thefinishstep. There is no green edge, and each $7$-cycle has at least five ports.}\par\nopagebreak[4]
		
		Suppose $U h V$ is green and $VW$ is black. We show that $W$ is an end of a green edge. Let $t$ be the neighbor of $V$ other than $h,W$. If all edges at $W$ are black, take $R=\{h,V,W\}$, $M=\{t\}$, and $S=\{V\}$. By girth and Step~\ref{lem:no-redgreen}, the boundary consists of four distinct cubic vertices. Thus $\delta(H)\ge2$, $s(H)=1$, and $b_8(H)=0$. Step~\ref{lem:red-matching} requires four ports for an unmarked $C_7$, and $|N(R)\setminus M|=3$, so $c_7(H)=0$. Hence $w(G)-w(H)+2s(H)=(13-3)+2=12=12|S|$.
		
		Suppose a red path $W p q X$ meets $W$. Let $b$ be the neighbor of $W$ other than $V,p$. By Steps~\ref{lem:no-redgreen} and~\ref{lem:red-matching}, $t,b$ and the two neighbors of $X$ other than $q$ are cubic. Take
		\[
		R=\{h,V,W,p,q,X\},\qquad M=\{t\},\qquad S=\{V,q\}.
		\]
		The set $S$ dominates $R\cup M$ and $w(R)=27$. First suppose $UX\notin E(G)$. If $N(V)\cap N(X)=\emptyset$, the four unmarked boundary vertices are distinct. They are $U,b$ and the two surviving neighbors of $X$. Hence $w(G)-w(H)=27-4=23$. The component of $t$ has minimum degree at least two and a single mark, giving $2s(H)-b_8(H)=2$. An unmarked $C_7$ cannot contain both surviving neighbors of $X$ as ports, since they share the deleted neighbor $X$. This leaves at most three ports, so $c_7(H)=0$. If $V,X$ have a common neighbor, it must be $t$ by girth. Then $d_H(t)=1$ and only three unmarked cubic vertices lose an edge. Thus $w(G)-w(H)=27-3=24$ and $c_7(H)=0$.
		
		If $UX\in E(G)$, let $Q=U h V W p q X U$. Keep $R,M,S$, and let $a,c$ be the neighbors of $U,X$ off $Q$, respectively. The boundary vertices $U,t,b,c$ are distinct by girth. In particular, $tX\notin E(G)$, since $U h V t X U$ would be a $5$-cycle. The vertex $U$ loses two edges, each of $b,c$ loses one, and $t$ is marked. Thus $w(G)-w(H)=27-(8-4)-2(5-4)=21$. The only leaf of $H$ is $U$, and its neighbor $a$ is still cubic.
		
		The all-black case at the start of this step, applied to $Ua$, implies that $a$ is incident with a red or green edge. The matching properties and Step~\ref{lem:no-redgreen} permit exactly one such edge. Let $z$ be the first internal vertex of its path. Then $z$ is unmarked and $d_H(z)=2$. Indeed, the vertex beyond $z$ cannot be a deleted root by the same restrictions, and each deleted internal vertex already has both neighbors on $Q$. The other black neighbor $y$ of $a$ survives, since an edge to $V,W$ or $X$ would give a cycle of length at most five. If $y=t$ or $y=c$, we obtain a $5$-cycle or a $4$-cycle, respectively. The equality $y=b$ is possible, giving a $6$-cycle and $d_H(y)=2$.
		
		The component containing $U$ has one leaf, whose cubic neighbor $a$ has an unmarked degree-two neighbor $z$. Thus \leafpart{lem:L3} applies, including when $y=b$. We apply~\eqref{eq1} to the other components. None can be an unmarked $C_7$ because only $U,b,c$ can be unmarked attachments. Adding the bounds gives $12\gd(H)\le w(H)-3$, so
		\[
		12\gd(G)\le24+12\gd(H)\le24+w(H)-3=w(G).
		\]
		
		If the green matching is nonempty, every black edge from one of its ends leads to another end, and Step~\ref{lem:no-redgreen} excludes red edges at its ends. Connectedness then makes the green matching perfect, with no red edges. Let $S$ consist of the internal vertices of the green paths. Their neighbors are the $2|S|$ roots, so $S$ dominates $G$. Hence $12\gd(G)\le12|S|\le4(2|S|)+5|S|=w(G)$, contrary to the choice of $G$.
		
		Only black and red edges remain. If a $7$-cycle has $k$ ports, its other $7-k$ vertices lie in pairs on red arcs. These arcs have distinct ends because the red edges form a matching. Thus $7-k$ is even and $7-k\le k$, so $k\in\{5,7\}$.
		
		\smallskip
		\noindent\refstepcounter{finishstep}\label{st:lastred}\emph{Step \thefinishstep. There is no red edge.}\par\nopagebreak[4]
		
		Suppose $VW\in E(G)$ joins distinct red paths $V r s Y$ and $W p q Z$. Let $t,u$ be the remaining black neighbors of $V,W$, respectively. Put $R=\{V,W,r,s,p,q\}$. Its four boundary vertices $Y,Z,t,u$ are distinct. The equality $Y=Z$ would contradict the red matching, a coincidence between a red end and a black neighbor would give a $4$- or $5$-cycle, and $t=u$ would give a triangle. Each boundary vertex is cubic and loses one edge. With $M=\emptyset$ and $S=\{r,p\}$, we have $w(G)-w(H)=28-4=24$. By Step~\ref{lem:no-green}, an unmarked $C_7$ needs at least five ports, so $c_7(H)=0$. This excludes $VW$.
		
		Now let $u p q v$ be a red path, with $N(u)\setminus\{p\}=\{u_1,u_2\}$ and $N(v)\setminus\{q\}=\{v_1,v_2\}$. These pairs are disjoint because a common vertex would give a $5$-cycle. All four vertices are cubic. The preceding paragraph excludes a red edge at any of them, so all their neighbors are cubic as well.
		
		Suppose, after relabeling if necessary, that $u_1v_1\in E(G)$. Take $R=\{u,v,u_1,v_1,p,q\}$, $M=\{u_2,v_2\}$, and $S=\{u,v\}$. The remaining neighbors of $u_1,v_1$ form $N(R)\setminus M$, and all four members of $N(R)$ are distinct by girth. Thus $w(G)-w(H)=26-2=24$ and $c_7(H)=0$.
		
		If no $u_iv_j$ is an edge, retain $S=\{u,v\}$ and take $R=\{u,v,u_1,u_2,p,q\}$ and $M=\{v_1,v_2\}$. Neither $u_1$ nor $u_2$ meets $v_1$ or $v_2$. The absence of a $4$-cycle excludes a second common neighbor of $u_1,u_2$. Thus there are four unmarked boundary vertices, each losing one edge. By Step~\ref{lem:no-green}, $5c_7(H)\le4$, so $c_7(H)=0$. Since $\delta(H)\ge2$ and $H$ has exactly two marks, $s(H)\ge1$ and $b_8(H)=0$. Hence $w(G)-w(H)+2s(H)\ge(26-4)+2=24=12|S|$. This gives~\eqref{eq2} and excludes the last red path.
		
		Thus $G$ is cubic. Let $n=|V(G)|$, choose $v\in V(G)$, and put $R=\{v\}\cup N(v)$. Leave $H=G-R$ unmarked. By girth, $G[R]$ is a star and the six edges leaving $R$ have distinct ends in $H$. Hence six vertices of $H$ have degree two and $n-10$ have degree three, giving
		\[
		w(H)=6\times 5+4(n-10)=4n-10.
		\]
		A $7$-cycle component of $H$ would require an edge to $R$ at each of its seven vertices. Only six edges leave $R$, and $H$ has no marks, so $c_7(H)=s(H)=b_8(H)=0$. By minimality, $12\gamma(H)\le w(H)$.
		
		Adding $v$ to a minimum dominating set of $H$ dominates $G$. Thus
		\[
		12\gamma(G)\le12+12\gamma(H)\le12+w(H)=4n+2.
		\]
		Since $3\gamma(G)$ is an integer, this gives $3\gamma(G)\le\lfloor n+\tfrac12\rfloor=n$. Therefore $\gamma(G)\le\lfloor n/3\rfloor$, contrary to $12\gd(G)>w(G)=4n$. This proves Proposition~\ref{prop1} and Theorem~\ref{thm:cubic}.
	\end{proof}
	
	\section*{Declarations}
	
	\paragraph*{Competing interests}
	The authors declare that they have no known competing financial interests or personal relationships that could have appeared to influence the work reported in this paper.
	
	\paragraph*{Data availability}
	No datasets were generated or analyzed during this study.
	
	\paragraph*{Use of generative AI}
	The authors used ChatGPT (OpenAI) to assist with language editing and to improve the presentation of the manuscript. The authors reviewed and revised the resulting text and take full responsibility for the content of the manuscript.

\end{document}